\documentclass[12pt,a4paper]{amsart}
\usepackage{amsmath,amssymb,amscd,amsthm,amsfonts}
\usepackage{epsfig}
\usepackage{hyperref}

\newtheorem{theorem}{Theorem}[section]

\newtheorem{corollary}{Corollary}[section]

\newtheorem{prop}{Proposition}[section]

\newtheorem{question}{Open problem}[section]

\newcommand{\tri}{\Delta}

\newcommand{\al}{\alpha}
\newcommand{\be}{\beta}
\newcommand{\ga}{\gamma}

\newcommand{\eps}{\varepsilon}
\newcommand{\R}{\mathbb{R}}
\newcommand{\A}{\mathcal{A}}

\newcommand{\diam}{\mathrm{diam}}
\newcommand{\conv}{\mathrm{conv\;}}

\newcommand{\aff}{\mathrm{aff\;}}
\newcommand{\rank}{\textrm{rank\;}}

\numberwithin{equation}{section}

\begin{document}

\title{Tverberg plus minus}

\author[B\'ar\'any]{Imre B\'ar\'any}
\thanks{}
\address{Alfr\'ed R\'enyi Institute of Mathematics \\ Hungarian Academy of Ciencies \\ H-1364 Budapest Pf. 127 Hungary} 
\address{Department of Mathematics \\
University College London \\
Gower Street, London, WC1E 6BT, UK} 
\email{barany.imre@renyi.mta.hu}
\author[Sober\'on]{Pablo Sober\'on}
\thanks{} 
\address{Mathematics Department, Northeastern University, 360 Huntington Ave., 
Boston, MA 02115, USA } 
\email{p.soberonbravo@northeastern.edu}

\keywords{Tverberg's theorem, sign conditions}
\subjclass[2000]{Primary 52A37, secondary 52B05}

\begin{abstract} We prove a Tverberg type theorem: Given a set $A \subset \R^d$ in general position with $|A|=(r-1)(d+1)+1$ and $k\in \{0,1,\ldots,r-1\}$, there is a partition of $A$ into $r$ sets $A_1,\ldots,A_r$ (where $|A_j|\le d+1$ for each $j$) with the following property. There is a unique $z \in \bigcap_{j=1}^r \aff A_j$ and it can be written as an affine combination of the element in $A_j$: $z=\sum_{x\in A_j}\al(x)x$ for every $j$ and exactly $k$ of the coefficients $\al(x)$ are negative. The case $k=0$ is Tverberg's classical theorem.

\end{abstract}

\maketitle

\section{Introduction and main result}\label{introd}

Assume $A=\{a_1,\ldots,a_n\}\subset \R^d$ where $n=(r-1)(d+1)+1$ and $r\ge 2$, $d\ge 1$ are integers. Suppose further that the coordinates of the $a_i$ (altogether $dn$ real numbers) are algebraically independent. A partition $\A =\{A_1,\ldots,A_r\}$ of $A$ is called {\sl proper} if $1\le |A_j|\le d+1$ for every $j \in [r]$. Here and in what follows $[r]$ stands for the set $\{1,\ldots,r\}$. We will show later (Proposition~\ref{cap}) that in this case the intersection of the affine hull of the $A_j$s is a single point $z$, that is, $\{z\}=\bigcap_{j=1}^r \aff A_j$. Equivalently, the following system of linear equations has a unique solution:
\begin{equation}\label{aff}
z=\sum_{x\in A_j}\al(x)x \quad \mbox{ and } \quad 1=\sum_{x\in A_j}\al(x) \quad \mbox{ for all }j \in [r].
\end{equation}

One form of Tverberg's classical theorem \cite{tv} puts extra conditions on the coefficients $\al(x)$ (consult \cite{Mat03} and the references therein for an introduction to the subject).

\begin{theorem}[Tverberg's theorem]\label{th:tverb} Under the above conditions there is a proper partition of $A$ into sets $A_1,\ldots,A_r$ such that $\al(x)\ge 0$ for all elements $x \in A$. In other words, $\{z\}=\bigcap_{j=1}^r\conv A_j$.
\end{theorem}

This means that the unique solution to (\ref{aff}) has $\al(x)>0$ for all $x \in A$. Can we require here that exactly one (or two or more) of the $\al(x)$ are negative? A partial answer comes from the following theorem, which is the main result of this paper.

\begin{theorem}\label{th:pm} Assume $k \in \{0,1,\ldots,r-1\}$. Under the conditions of Theorem~\ref{th:tverb} there is a (proper) partition of $A$ into $r$ parts so that in the unique solution to (\ref{aff}) $\al(x) <0$ for exactly $k$ elements $x \in A$.
\end{theorem}

Of course the same holds for any set $A$ of $n$ points in $\R^d$, we only have to relax the condition $\al(x) <0$ to $\al(x)\le 0$ for $k$ elements $x \in A$ and $\al(x)\ge 0$ for the rest. Actually, the general position condition (on $A$) is used in order to avoid cases when $\al(x)=0$ for some $x \in A$.

It is not clear for what other values of $k$, $k\in [n]$, the theorem holds.  Since the sum of the coefficients for each $A_j$ is one, at least one is positive. This implies the upper bound $k\le n-r$.

The case $d=1$ is very simple. Then $n=2r-1$ and there is no $r$-partition with $r$ or more negative coefficients, so the trivial bound $k\le n-r=r-1$ is tight. In the case $d=2,\;r=3$ and $n=7$ Theorem~\ref{th:pm} gives a suitable partition for $k=0,1,2$. A careful case analysis shows that the statement holds for $k=3$ as well, and an extensive  computer aided search did not find any example where it fails to hold for $k=4$.

The case of $r=2$, that is, Radon (plus minus) partitions can be checked directly. Then $|A|=d+2$ and the outcome is that for any $k\in \{0,1,\ldots,\lfloor \frac {d+2}2\rfloor\}$ there is a partition with exactly $k$ negative $\al(x)$. Further, there are examples showing that this does not hold for $k>\lfloor \frac {d+2}2\rfloor$.  In this case everything is governed by the unique affine dependence of the vectors in $A$, just as in the proof of Radon's theorem. We omit the details.

We will see in Corollary~\ref{miracle} in Section 3 that, for a strange reason, if both $d$ and $r$ are even, then Theorem~\ref{th:pm} holds with $k=\frac 12 [(r-1)(d+1)+1]$ as well. This makes us wonder if Theorem~\ref{th:pm} holds for all integers $k\le \frac 12 [(r-1)(d+1)+1]$.

We are going to prove Theorem~\ref{th:pm} in a stronger form: to some extent we can prescribe the subset of $A$ where the coefficients in (\ref{aff}) are negative.

\begin{theorem}\label{th:tverbM} Under the conditions of Theorem~\ref{th:tverb} let $M\subset A$ be a set of size at most $r-1$ such that $\conv M \cap \conv (A\setminus M)=\emptyset$. Then there is a partition  $\A =\{A_1,\ldots,A_r\}$ of $A$ such that in (\ref{aff}) $\al(x)<0$ if and only if $x \in M$.
\end{theorem}

We prove this theorem in Section 3 where we state a slightly stronger result whose proof is in Section 6. Examples showing the necessity of the condition on $M$ are given in Section 2. In Section \ref{section-coloured} we discuss coloured variations of Theorem \ref{th:pm}. In Section 5 we prove the following fact.

\begin{prop}\label{cap} Assume $A=\{a_1,\ldots,a_n\}\subset \R^d$, the coordinates of the $a_i$ are algebraically independent and $r\ge 2$, $d\ge 1$ are integers. If the partition $\A =\{A_1,\ldots,A_r\}$ of $A$ is proper and  $n=(r-1)(d+1)+1$, then $\bigcap_{j=1}^r \aff A_j$ is a single point. If $n\le (r-1)(d+1)$, then $\bigcap_{j=1}^r \aff A_j=\emptyset$.
\end{prop}

The last statement holds even if the partition is not proper. The first part must be known, see for instance \cite{PS} or \cite{dv} for similar statements. The second part is proved in \cite{tv}. We give a simple proof in Section \ref{section-proposition}.

\section{The condition on $M$}\label{prep}

The condition on $M$ in Theorem \ref{th:tverbM} simply says that $M$ and $A \setminus M$ can be separated by a hyperplane.

{\bf Example 1.} We give an example showing the necessity of this condition. Let $V=\{v_1,\ldots,v_{d+1}\}$ be the set of vertices of a regular simplex $\tri$, and let $c$ be the centre of $\tri$ and write $F_h$ for its facet opposite to $v_h$. For every $h \in [d+1]$ let $U_h\subset v_h+\eps B$ be an $(r-1)$-element set. Here $\eps>0$ is small and $B$ is the Euclidean unit ball in $\R^d$ centred at the origin. Define $A=\{c\}\cup \bigcup_{h=1}^{d+1}U_h$. We assume $A$ is in general position which can be clearly reached by choosing the sets $U_h$ suitably. Set $M=\{c\}$ so the separation condition fails. We {\bf claim} that there is no proper $r$-partition of $A$ such that in (\ref{aff}) only $\al(c)$ is negative, $\al(x)>0$ for  all other $x \in A$.

Assume the contrary and let  $A_1,\ldots,A_r$  be a proper partition with $z \in \bigcap_{j=1}^r \aff A_j$ so that $\al(c) < 0$ and $\al(x)> 0$ for all other $x \in A$ in (\ref{aff}).

\begin{figure}[h]
\centerline{\includegraphics[scale=1]{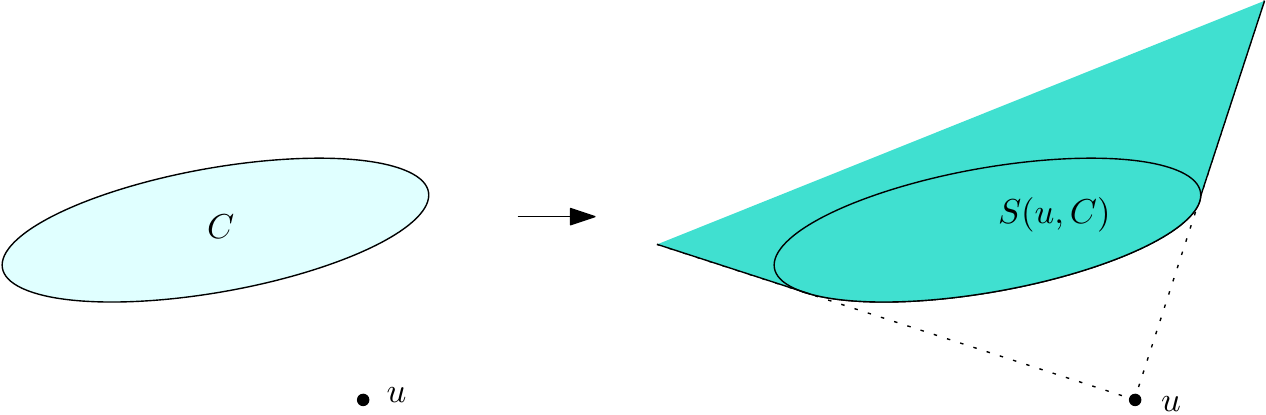}}
\caption{Construction of $S(u,C)$.  Notice that $C \subset S(u,C)$.}\label{figure-shadows}
\end{figure}

Given a convex compact set $C$ in $\R^d$ and a point $u \in \R^d \setminus C$ we let $S(u,C)$ denote the {\sl shadow} of $C$ from $u$ which is the set of point $\{tu+(1-t)c: c\in C,\; t\le 0\}$, see Figure \ref{figure-shadows}.

\begin{figure}[h]
\centerline{\includegraphics[scale=1]{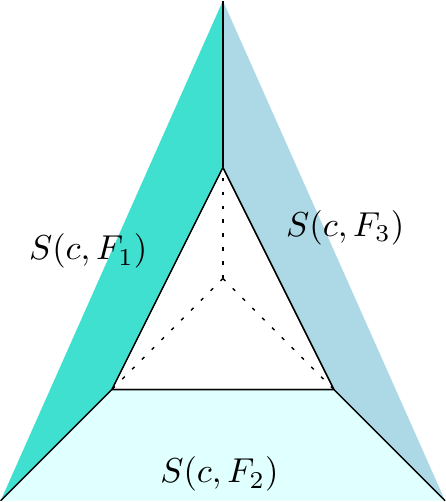}}
\caption{The shadows of the sides of a triangle from its centre do not intersect.}
\label{figure-triangle}
\end{figure}

For every $h \in [d+1]$ there is a $j=j(h) \in [r]$ such that $A_{j(h)}$ and $U_h$ are disjoint, simply because each $U_h$ has $r-1$ elements and the number of sets $A_j$ is $r$. It follows that $A_{j(h)}\subset F_h + \eps B \subset S(c, F_h+\eps B)$  if $c \not\in A_{j(h)}$, and then $z \in \conv A_{j(h)} \subset S(c, F_h+\eps B)$. If $c\in A_{j(h)}$, then  in the equation $z=\sum_{x\in A_{j(h)}}\al(x)x$ only the coefficient $\al(c)$ is negative, so $z \in S(c, F_h+\eps B)$. Therefore

\[
z \in \bigcap_{h=1}^{d+1} S(c,F_h+\eps B). 
\]
However, $ \bigcap_{h=1}^{d+1} S(c,F_h+\eps B)=\emptyset$ as long as $\eps < \frac{\diam \tri}{2d}$.  See Figure \ref{figure-triangle} for an illustration.  This follows from the fact that the shadows $S(c,F_h+\eps B)$ for $h \in [d+1]$ are convex and their union covers the boundary of $\tri$.  If they had a point in common, then their union would cover $\tri$.  Since none of them contains $c$, this is impossible.  This gives us the contradiction we were seeking. 

\bigskip
The proof above does not use the fact that $c$ is the only point near the centre of the simplex.  If we consider $U_0 \subset c + \eps B$ any set and declare $U_0 = M$, the same arguments as above show that there is no partition with the desired properties.  Therefore, the condition $\conv M \cap \conv (A \setminus M)= \emptyset$ cannot be removed even if we allow $A$ to have more than $(r-1)(d+1)+1$ points.

\bigskip
{\bf Example 2.} This example shows a construction where $M$ satisfies the separation condition, $|M|=r$ and the conclusion of the theorem fails. We work with the same simplex $\tri$ and $U_h\subset v_h+\eps B$ is the same $(r-1)$-element set as before for $h\in \{2,\ldots,d+1\}$ but for $h=1$ it is an $r$-element set in $v_1+\eps B$. This time $A=\bigcup_{h=1}^{d+1}U_h$ and $M=U_1$. We assume of course that $A$ is in general position. Now $M$ is separated from $A \setminus M$ and has exactly $r$ elements.

We {\bf claim} that $A$ has no partition into $r$ parts with the required properties. 

{\bf Proof.} Assume the contrary and let  $A_1,\ldots,A_r$  be a proper partition with $z \in \bigcap_{j=1}^r \aff A_j$ such that in (\ref{aff}) $\al(x) < 0$ if $x \in M$ an $\al(x)>0$ if $x \notin M$. 
For $h \in \{2,\dots,d+1\}$ let $G_h$ be the convex hull of $V\setminus\{v_1,v_h\}$; this is a $(d-2)$-face of $\tri$. Set $\beta=\sum_{x \in A_j\cap U_1}\al(x)$, and note that $\beta>0$ if $A_j \cap U_1$ is nonempty. Define $u_j=\frac 1{\beta} \sum_{x \in A_j\cap U_1}\al(x)x$ if $A_j \cap U_1$ is nonempty, and $u_j=v_1$ otherwise. 

Note again that for each $h\in \{2,\ldots,d+1\}$ there is a $j(h)\in [r]$ such that $U_h$ and $A_{j(h)}$ are disjoint. Then $A_{j(h)}\setminus U_1 \subset G_h+\eps B$. This and the sign condition in (\ref{aff}) imply that for every $h>1$ with $A_{j(h)}\cap U_1\ne \emptyset$, 
\[
z \in S(u_{j(h)},G_h+\eps B). 
\]
This holds even if $A_{j(h)}\cap U_1= \emptyset$  since then $u_{j(h)}=v_1$ and $z \in \conv A_{j(h)}\subset G_h+\eps B \subset S(v_1,G_h+\eps B)$. Thus 
\[
z \in \bigcap_{h=2}^{d+1} S(u_{j(h)},G_h+\eps B). 
\]
But again, the shadows on the right hand side have no point in common, as one can check easily.\qed

\section{Proof of Theorem~\ref{th:tverbM}}\label{section-proof}

We are going to use the colourful Carath\'eodory theorem \cite{bar82}. It says that given sets $S_1,\ldots,S_{n+1} \subset \R^n$ with the condition that $0 \in \bigcap_{i=1}^{n+1} \conv S_i$, there is a transversal, that is, a choice $s_i \in S_i$ for every $i \in [n+1]$, such that $0 \in \conv \{s_1,\ldots,s_{n+1}\}$.

{\bf Proof} of Theorem~\ref{th:tverbM}. We use a modification of Sarkaria's argument~\cite{sark}, in the form given by~\cite{BarOnn}. It starts with an artificial tool: let $v_1,\ldots,v_r$ be the vertices of a regular simplex in $\R^{r-1}$ centred at the origin. The important property here is that, apart from scalar multiples, their unique linear dependence is $v_1+\ldots+v_r=0$.

Assume $A=\{a_1,\ldots,a_n\}$ where $n=(r-1)(d+1)+1$. Recall that $M \subset A$, $|M|=k < r$, and $M$ and $A \setminus M$ are separated by a hyperplane. Define 
\[
b_i=(a_i,1) \mbox{ if } a_i \notin M \mbox{ and } b_i=(-a_i,-1) \mbox{ if } a_i \in M.
\]
where $(a_i,1)\in \R^{d+1}$ is vector $a_i$ appended with an $(d+1)$-th coordinate equal to one, and similarly for $(-a_i,-1)$. For $i \in [n]$ we set
\[
S_i=\{v_1\otimes b_i,v_2\otimes b_i,\ldots,v_r\otimes b_i\}.
\]
Here $v_j\otimes b_i$ is the usual tensor product, which is the same as the matrix product of the $(r-1)$-dimensional column vector $v_j$ and the $(d+1)$-dimensional row vector $b_i^T$: $v_j b_i^T$, where we consider our vectors as vertical matrices. So this product is an $(r-1)\times (d+1)$ matrix, or equivalently a vector in $\R^{n-1}$. Observe that $0 \in \conv S_i$ for every $i$, so the colourful Carath\'eodory theorem applies and gives a transversal $s_1,\ldots,s_n$ whose convex hull contains the origin, that is, there are non-negative coefficients $\be_1,\ldots,\be_n$ whose sum is 1 such that $\sum_{i=1}^n\be_is_i=0$. Here each $s_i$ is of the form $v_j\otimes b_i$ for a unique $j=j(i) \in [r]$. We define $A_j=\{a_i\in A: j(i)=j\}$ for all $j \in [r]$. The sets $A_1,\ldots,A_r$ form an $r$-partition of $A$. With the new notation
\begin{eqnarray*}\label{eq:crux}
0&=&\sum_{i=1}^n\be_is_i=\sum_{i=1}^n\be_i v_{j(i)}\otimes b_i\\
  &=&\sum_{j=1}^r\sum_{a_i \in A_j}\be_iv_j\otimes b_i=\sum_{j=1}^rv_j\otimes \left(\sum_{a_i \in A_j}\be_i b_i\right).
\end{eqnarray*}

 Define now $\al_i=-\be_i$ if $a_i \in M$ and  $\al_i=\be_i$ otherwise. The last equation becomes
\begin{equation}\label{eq:crux}
0=\sum_{j=1}^rv_j\otimes \left(\sum_{a_i \in A_j}\al_i (a_i,1)\right).
\end{equation}
There is a vector $u \in \R^{r-1}$, orthogonal to $v_3,v_4,\ldots,v_r$ with $\langle u, v_1\rangle =1$, where $\langle \cdot , \cdot \rangle$ denotes the dot product. The condition  $v_1+\ldots+v_r=0$ implies that $\langle u, v_2\rangle =-1$. As (\ref{eq:crux}) is a matrix equation, multiplying it from the left by the $(r-1)$-dimensional row vector $u^T$ gives $\sum_{a_i \in A_1}\al_i (a_i,1)=\sum_{a_i \in A_2}\al_i (a_i,1)$.  By symmetry we have
\begin{equation}\label{equal}
z:=\sum_{a_i \in A_1}\al_i (a_i,1)=\sum_{a_i \in A_2}\al_i (a_i,1)=\ldots =\sum_{a_i \in A_r}\al_i (a_i,1).
\end{equation}
There are two cases to be considered.

\medskip
{\bf Case 1:} when $A_j=\emptyset$ for some $j \in [r]$. Then $z=0$ and some $A_h$, say $A_1$, is nonempty and not all coefficients $\al_i$ with $a_i \in A_1$ are zero. Thus $\sum_{a_i \in A_1}\al_i (a_i,1)=0$.  Then  $\al_i\le 0$ for all $a_i \in A_1\cap M$ and  $\al_i\ge 0$ for all  $a_i \in A_1\setminus M$. Setting
\[
\ga :=\sum_{a_i \in A_1\cap M}\al_i=\sum_{a_i \in A_1\setminus M}-\al_i,
\]
it follows that $\ga >0$. Consequently $\conv (A_1\cap M)$ and $\conv (A_1\setminus M)$ have a point in common, namely 
\[\frac 1{\ga}\sum_{a_i \in A_1\cap M}\al_i a_i=\frac 1{\ga}\sum_{a_i \in A_1\setminus M}(-\al_i) a_i,\] contradicting the separation assumption.

\medskip
{\bf Case 2:} when $A_j$ is nonempty for all $j \in [r]$ . Reading the last coordinate of (\ref{equal}) gives that
\[
\ga:= \sum_{a_i \in A_1}\al_i=\sum_{a_i \in A_2}\al_i =\ldots =\sum_{a_i \in A_r}\al_i.
\]
Since $|M|<r$, there is a $j \in [r]$ such that $\al_i>0$ for all $a_i \in A_j$, implying that $\ga>0$. Then the point $\frac{1}{\ga}z$ is in the affine hull of every $A_j$. The construction guarantees that $\al_i<0$ if and only if $a_i \in M$. \qed

\bigskip
Actually, this proof gives a stronger statement. In Case 2, the positivity of $\ga$ is guaranteed by the condition $|M|=k<r$. Not assuming $k<r$, $\ga$ can be negative or zero. When $\ga<0$ equation (\ref{equal}) implies again that $\frac{1}{\ga}z$ is in the affine hull of every $A_j$, but this time $\al(x)>0$ exactly when $x \notin M$. We will exclude the case $\ga=0$ using the general position condition. The proof of this is given in Section~\ref{nonzero} because it uses the content of Section~\ref{section-proposition}. So we have the following result.

\begin{theorem}\label{th:stverbM} Under the conditions of Theorem~\ref{th:tverb} let $M$ be a subset $A$ such that $\conv M \cap \conv (A\setminus M)=\emptyset$. Then there is a partition  $\A =\{A_1,\ldots,A_r\}$ of $A$ such that in (\ref{aff}) either 
\begin{itemize}
	\item $\al(x)<0$ if and only if $x \in M$, or 
         \item $\al(x)<0$ if and only if $x \notin M$.
\end{itemize}
\end{theorem}

The second example in Section~\ref{prep} shows that in some cases only the second alternative holds.

\begin{corollary} \label{miracle}Assume $r,d$ are both even and positive integers, $A\subset \R^d$ is in general position, $|A|=(r-1)(d+1)+1$, and $k=\frac 12 [(r-1)(d+1)+1]$. Then $A$ has a proper $r$-partition such that in (\ref{aff}) exactly $k$ of the coefficients $\al(x)$ are negative. 
\end{corollary}

The {\bf proof} is easy. Under the above conditions there is a subset $M$ of $A$ of size $k$ that is separated from $A\setminus M$. According to Theorem~\ref{th:stverbM}, $A$ has an $r$-partition such that in (\ref{aff}) either $\al(x)<0$ if and only if $x \in M$, or $\al(x)<0$ if and only if $x \in A\setminus M$. In both cases, exactly $k$ coefficients are negative.\qed

\medskip
{\bf Remark.} The same result can be proved using Tverberg's original method of {\sl moving the points}.  The main idea is to start with a set of points which have a partition with the required conditions.  Then, as one moves one point continuously, if the partition stops working, one can show that points may be swapped in the partition in order to still satisfy the conclusion of the theorem.  The proof given above is shorter and simpler.

\section{Colourful Tverberg plus minus}\label{section-coloured}

Once a Tverberg type theorem with conditions on the signs of coefficients of the affine combinations has been established, it becomes natural to try to extend it to the coloured versions of Tverberg's theorem, as in \cite{BL92}.

Given disjoint sets $F_1, \ldots, F_n$ of $r$ points each in $\R^d$, considered as colour classes, we say that $A_1, \ldots, A_r$ is a colourful partition of them if $|F_i \cap A_j| = 1$ for all $i\in [n]$, $j \in  [r]$.  In such a case, we can denote the points by $x_{i,j} = F_i \cap A_j$.  The coloured Tverberg theorem is concerned about the existence of colourful partitions for which the convex hulls of the sets $A_j$ intersect. In other words, we seek a colourful partition and a point $z\in \R^d$ for which there is a solution to the equations
\begin{eqnarray}
z & = & \sum_{i=1}^n \alpha(x_{i,j})x_{i,j} \mbox{ for all }j \in [r] \label{condition-1}\\
&&\mbox{subject to } 1 = \sum_{i=1}^n \alpha (x_{i,j}) \mbox{ for all } j \in [r], \ \mbox{and} \label{condition-2}\\
& & \alpha({x_{i,j}}) \ge 0 \mbox{ for all } i \in [n], j \in [r]. \label{condition-coloured}
\end{eqnarray}

The question then becomes, given $M \subset [n]$, find a solution where we exchange (\ref{condition-coloured}) for
\begin{eqnarray}
\alpha({x_{i,j}})  \le 0 & & \mbox{ for all } i \in M, j \in [r], \mbox{and} \\
\alpha({x_{i,j}}) \ge 0 & & \mbox{ for all } i \in [n]\setminus M, j \in [r].
\end{eqnarray}

In other words, we aim to prescribe negative coefficients, but we also require that the same restrictions hold accross the colour classes.  We obtain a partial result in this direction.

\begin{theorem}\label{theorem-colouredpm}
Let $n=(r-1)d+1$ and $F_1, \ldots, F_n$ be disjoint subsets of $\R^d$ whose union is algebraically independent, each of cardinality $r$ and $M \subset [n]$.  Then, there is a colourful partition of $F_1, \ldots, F_n$ into $r$ sets and solutions to equations (\ref{condition-1}) and (\ref{condition-2}) such that either
\begin{itemize}
	\item $\alpha(x_{i,j}) > 0$ for all $i \in M$ and $\alpha(x_{i,j})< 0$ for all $i \in [n]\setminus M$, or
	\item $\alpha(x_{i,j}) < 0$ for all $i \in M$ and $\alpha(x_{i,j})> 0$ for all $i \in [n]\setminus M$.
\end{itemize}
Moreover, the affine combinations use the same coefficients for the colour classes.  In other words, for all $i \in [n]$ and $j, j' \in [r]$, $\alpha(x_{i,j}) = \alpha(x_{i,j'})$.
\end{theorem}

{\bf Proof.}
We use the main result of \cite{Sob}.  It says that for $n = (r-1)d+1$ and the sets $F_1, \ldots, F_n$, there are solutions to equations (\ref{condition-1}), (\ref{condition-2}) and (\ref{condition-coloured}) where $\alpha(x_{i,j}) = \alpha (x_{i,j'})$ for all $i \in [n], j \in [r], j' \in [r]$.  Then, we apply this result to the sets \[G_i = \begin{cases}
F_i & \mbox{if } i \in M \\
-F_i & \mbox{otherwise.}
\end{cases}\]  Let $B_1, \ldots, B_r$ be the colourful partition we obtain of $G_1, \ldots, G_n$, with $y_{i,j} = G_i \cap B_j$ for all $i,j$.  We denote by $\beta (y_{i,j})$ the coefficients we obtain satisfying equations (\ref{condition-1}), (\ref{condition-2}) and (\ref{condition-coloured}).  We rename them as $\beta (y_{i,j}) = \beta_i$, since they do not depend on $j$.  Let $x_{i,j} = \pm y_{i,j}$ and $\alpha_i = \pm \beta_i$ where the sign is positive (negative) if $i \in M$ ($i \notin M$), respectively.  Let $\ga = \sum_{i=1}^n \alpha_i$.  Let us see what happens if $\gamma \neq 0$.  By construction, if we consider $A_1, \ldots, A_r$ the partition induced by the points $x_{i,j}$ and $\alpha(x_{i,j}) = \alpha_i / \ga$ for all $i,j$, they satisfy all the requirements for the conclusion of the theorem.  The two cases in Theorem \ref{theorem-colouredpm} correspond to the possibilities for the sign of $\ga$.

We have to verify that the general condition assumption we have on $F_1, \ldots, F_n$ implies $\ga \neq 0$. This part of the proof is technical, and it relies on the modification of Sarkaria's trick from \cite{Sob}.  Given a set $F= \{z_1, \ldots, z_r\} \subset \R^d$, a permutation $\sigma: [r] \to [r]$ and $v_1, \ldots, v_r \in \R^{r-1}$ as in Section \ref{section-proof}, we can define
\[
F \otimes \sigma = \sum_{j=1}^r z_j \otimes v_{\sigma(j)} \in \R^{n-1} ,\qquad S(F) = \{F \otimes \sigma: \sigma \ \mbox{is a permutation}\}
\]
The existence of $\beta_1, \ldots, \beta_n$ follows from applying the colourful Carath\'eodory theorem to the sets $S(G_1), \ldots, S(G_n)$ in $\R^{n-1}$.  However, if the original set of points $\cup_{i=1}^n F_i$ is algebraically independent, then no transversal to $S(F_1), \ldots, S(F_n)$ would have a non-trivial affine dependence in $\R^{n-1}$, so $\gamma \neq 0$, as required.
\qed

As mentioned, Theorem \ref{theorem-colouredpm} has equal coefficients accross the colour classes.  Removing this condition leads to the following problem.

\begin{question}
If we remove the equal coefficients condition, does Theorem \ref{theorem-colouredpm} hold with $n=d+1$?
\end{question}

The answer is affirmative with $r=2$.  If $M = \emptyset$ this is the main conjecture from \cite{BL92}.

\section{Proof of Proposition~\ref{cap}}\label{section-proposition}
We write the equation (\ref{aff}) in matrix form $M\al =b$. The $(n+d)\times(n+d)$ matrix $M$ is made up of blocks. The block corresponding to $A_j$ is a $d\times |A_j|$ matrix $N_j$ whose columns are the vectors in $A_j$. The row immediatley below block $N_j$ has a 1 in each column containing a vector from $A_j$ and zeroes everywhere else. There are $r$ further blocks, each one is $-I_d$, the negative $d \times d$ identity matrix. They are in the last $d$ columns of $M$, with a row of zeroes between them. These submatrices are arranged in $M$ as shown on Table 1. All other entries of $M$ are zeroes. The $i$th column of $M$ corresponds to the vector $a_i$. Note that $M=M(\A)$ depends on $A$ and on the partition $\A=\{A_1,\ldots,A_r\}$ as well.

\begin{table}[t]\label{fig:1}
\begin{center}
\begin{tabular}{|c|c|c|c|c|}
  \hline
   & & & & \\
   $N_1$ & & & & \hspace*{.1in} $-I_d$
                           \hspace*{.1in}  \\
   & & & & \\
   \hline
     1 1 $\dots$ 1 & & & & \\
  \hline
   & & & & \\
  &  $N_2$ & & & $-I_d$ \\
   & & & & \\
   \hline
    & 1 1 $\dots$ 1 & & & \\
  \hline
  & & $\ddots$ & & \\
  \hline
   & & & & \\
   & & & $N_r$ & $-I_d$ \\
   & & & & \\
  \hline
     & & & 1 1 $\dots$ 1  & \\
     \hline
\end{tabular}
\bigskip
\end{center}
\caption{The matrix $M$, the empty regions indicate zeros}
\end{table}

The variables are $\al =(\al_1,\ldots,\al_n,z_1,\ldots,z_d)^T \in \R^{n+d}$ and the right hand side vector is $b\in \R^{n+d}$ that has coordinate zero everywhere except in positions $d+1,2(d+1),\ldots,r(d+1)$ where it has one. The original system (\ref{aff}) is the same as
\begin{equation}\label{matrix}
M\al =b.
\end{equation}

As we have seen, $\bigcap_{j=1}^r \aff A_j$ is a single point if and only if the linear system (\ref{aff}), or what is the same, the equation~(\ref{matrix}) has a unique solution which happens if and only if $\det M \ne 0$. Here $\det M$ is a polynomial with integral coefficients in the coordinates of the $a_i$. If this polynomial is zero at some algebraically independent points $a_1,\ldots,a_n$, then it is identically zero. So it suffices to show one example where it is non-zero or, what is the same, one example where $\bigcap_1^r \aff A_j$ is a single point.

The example is simple.  Suppose $|A_j|=d+1-m_j$ for all $j\in [r]$ and $m_1\ge m_2\ge \ldots \ge m_r$. As $\A$ is a proper partition, $0\le m_j\le d$. Let $H_j$ be the subspace of $\R^d$ defined by equations $x_i=0$ for $i=\sum_{h=1}^{j-1}m_h+1,\ldots,\sum_{h=1}^jm_h$. Since $n=(r-1)(d+1)+1$, $\sum_1^r m_j=d$, implying that $\bigcap_1^r H_j$ is a single point, namely the origin. For each $p\in [r]$ choose $|A_j|$ affinely independent points in $H_j$. Their affine hull is exactly $H_j$, finishing the proof of the first part.

For the second part we can assume that $A_j$ is nonempty for all $j$, and also that $|A_j|\le d+1$ as otherwise one can delete some elements of $A_j$ while keeping its affine hull the same. We suppose further that $n=(r-1)(d+1)$ by adding extra (and algebraically independent) points to some suitable $A_j$s. Then $\bigcap_{j=1}^r \aff A_j\ne \emptyset$ if and only if the corresponding linear system (\ref{matrix}) has a solution. Now $M$ is an $(n+1)\times n$ matrix. Adding $b$ to $M$ as a last column we get a matrix that we denote by $M^*$. The system (\ref{matrix}) has a solution if and only if $\rank M=\rank M^*$. The previous argument shows that $\rank M=n-1$ and so we have that, as a polynomial, $\det M^*$ is identically zero. Again it suffices to give a single example where $\bigcap \aff A_j=\emptyset$. We use the same example as before except that this time $\sum_{j=1}^r m_j=d+1$ so we can add the equation $\sum_{i=1}^dx_i=1$ to the ones defining $H_1$ if $m_1<d$ and then $\bigcap H_j=\emptyset$, indeed. If $m_1=d$ then $H_1=0$ and $m_2=1$ and we define $H_2$ by the single equation $x_1+x_2=1$, and again $\bigcap H_j=\emptyset$. The sets $A_j$ are constructed the same way as above.\qed

\section{Proof of Theorem~\ref{th:stverbM}}~\label{nonzero}

{\bf Proof.} As we have seen we only have to show that $\ga\ne 0$. Assume $\ga=0$. This happens if and only if the homogeneous version of equation~(\ref{matrix}), that is
\begin{equation}\label{matrix0}
A\al=0
\end{equation}
has a nontrivial solution, which happens again if and only if $\det M=0$. This is impossible if the partition is proper (as we have seen in the previous section). Note that $z\ne 0$ and no $A_j$ is the emptyset, this follows from Case 1 of the proof of Theorem~\ref{th:tverbM}. So assume the partition is not proper. Then $A_j$ has more than $d+1$ elements for some $j$. Assume that $|A_1|>d+1$, say. This means that
\[
\sum_{x \in A_1}\al(x)(x,1)=(z,0).
\]
The vectors  $(x,1)$, $x \in A_1$ are affinely dependent, implying that there is a non-trivial affine dependence $\sum_{x \in A_1}\be(x)(x,1)=(0,0)$. Then for all $t \in \R$
\[
\sum_{x \in A_1}(\al(x)+t\be(x))(x,1)=(z,0),
\]
We choose here $t=t_0$ so that $\al(x)+t_0\be(x)=0$ for some $x=x_0\in A_1$. Set $\al'(x)=\al(x)+t_0\be(x)$ when $x \in A_1$ and $\al'(x)=\al(x)$ otherwise.

We change now the partition $A_1,\ldots,A_r$ to another one $A_1'\ldots,A_r '$ as follows. Set  $A_1'=A_1\setminus \{x_0\}$ and choose  some $A_j$ with $|A_j|\le d$ and set $A_j'=A_j\cup\{x_0\}$. All other $A_h$ remain the same. Let $M'$ be the corresponding matrix. 

We claim now that $\det M'=0$. The linear system (\ref{matrix0}) has a nontrivial solution, namely $\al(x)=\al'(x)$ with $z\ne 0$ unchanged. So indeed, $\det M'= 0$. 

Repeating this step finitely many times gives a proper partition such that (\ref{matrix0}) has a nontrivial solution. But the previous section shows that for a proper partition, (\ref{matrix0}) has no non-trivial solution.\qed

\bigskip
{\bf Acknowledgements.}  This material is partly based upon work supported by the National Science Foundation under Grant No. DMS-1440140 while the first author was in residence at the Mathematical Sciences Research Institute in Berkeley, California, during the Fall 2017 semester. Support from
Hungarian National Research Grants no K111827 and K116769 is acknowledged. We are also indebted to Attila P\'or and Manfred Scheucher for useful discussions, and to an anonymous referee for careful reading and valuable comments.

\bigskip




\end{document}